\title{Dense graphs with a large triangle cover have a large triangle packing}

\author{
Raphael Yuster\\
Department of Mathematics\\
University of Haifa\\
Haifa 31905, Israel\\
E--mail: raphy@math.haifa.ac.il
}

\date{}

\documentclass [letterpaper,11pt]{article}

\newtheorem{theo}{Theorem}[section]
\newtheorem{coro}[theo]{Corollary}
\newtheorem{lemma}[theo]{Lemma}
\newtheorem{conj}[theo]{Conjecture}

\newcommand{\qed}{\hspace*{\fill} \rule{7pt}{7pt}}

\begin{document}
\maketitle

\begin{abstract}

It is well known that a graph with $m$ edges can be made triangle-free by removing (slightly 
less than) $m/2$ edges. On the other hand, there are many classes of graphs
which are hard to make triangle-free in the sense that it is {\em necessary}
to remove roughly $m/2$ edges in order to eliminate all triangles.

We prove that dense graphs that are hard to make triangle-free, have a large packing of pairwise edge-disjoint triangles.
In particular, they have more than $m(1/4+c\beta^2)$ pairwise edge-disjoint triangles where $\beta$ is the density of the graph
and $c$ is an absolute constant. This improves upon a previous $m(1/4-o(1))$ bound which follows from the
asymptotic validity of Tuza's conjecture for dense graphs.
We conjecture that such graphs have an asymptotically optimal triangle packing of size $m(1/3-o(1))$.

We extend our result from triangles to larger cliques and odd cycles.

\end{abstract}

\section{Introduction}

All graphs in this paper are finite, undirected, and simple.
A {\em triangle cover} in a graph is a set of edges meeting all triangles. In other words, the removal of a triangle cover
results in a triangle-free graph. Dually, a {\em triangle packing} in a graph is a set of pairwise edge-disjoint triangles.
We denote by $\tau_t(G)$ the minimum size of a triangle cover and by $\nu_t(G)$ the maximum size of a triangle packing of a graph $G$.
It is easily observed that:
$$
\nu_t(G) \le \tau_t(G) \le 3\nu_t(G)\;.
$$
The first inequality follows from the fact that one must delete at least one edge from each triangle in a triangle packing in order to
obtain a triangle-free graph. The second inequality follows from the fact that deleting all edges of all triangles in a maximum triangle packing
results in a triangle-free graph. A long standing conjecture of Tuza \cite{Tu81} states that this second inequality is not optimal.
\begin{conj}[Tuza \cite{Tu81}]
\label{con:tuza}
$\tau_t(G) \le 2\nu_t(G)$.
\end{conj}
This conjecture, if true, is best possible as can be seen by taking, say, $G=K_4$ or $G=K_5$.
The best upper bound for $\tau_t(G)$ is due to Haxell \cite{Ha99} who prove that $\tau_t(G) \le (3-\frac{3}{23})\nu_t(G)$.

However, there is an important setting, where, asymptotically, Tuza's conjecture  holds.
This is the {\em dense graph} setting. To state this result we need first to consider the fractional relaxations of $\tau_t(G)$ and $\nu_t(G)$.
A {\em fractional triangle cover} assigns nonnegative weights to the edges so that the resulting weight of each triangle
(being the sum of the weights of its edges) is at least $1$. Dually, a {\em fractional triangle packing} assigns nonnegative weights
to the triangles so that the resulting weight
of each edge (being the sum of the weights of the triangles it meets) is at most $1$.
The goal is thus to minimize the sum of the weights of a fractional triangle cover and to maximize the sum of the weights
of a fractional triangle packing.
Let, therefore, $\tau_t^*(G)$ and $\nu_t^*(G)$ be the fractional relaxations of $\tau_t(G)$ and $\nu_t(G)$ respectively.
By linear programming duality we have $\tau_t^*(G)=\nu_t^*(G)$. Krivelevich \cite{Kr95} proved that Tuza's conjecture holds in a mixed fractional-integral setting. Namely, he proved:
\begin{theo}[Krivelevich \cite{Kr95}]
\label{t:kriv}
For any graph $G$ we have $\tau_t(G) \le 2\nu_t^*(G)$ and $\tau_t^*(G) \le 2\nu_t(G)$.
\end{theo}
The inequality $\tau_t^*(G) \le 2\nu_t(G)$ is tight (e.g. $K_4$) and the inequality $\tau_t(G) \le 2\nu_t^*(G)$ is known to be asymptotically tight.
A few years later, Haxell and R\"odl \cite{HaRo01} (see also \cite{Yu05}) proved that $|\nu_t(G)-\nu_t^*(G)| = o(n^2)$ for $n$-vertex graphs $G$.
In other words, in graphs that contain a quadratic number of pairwise edge-disjoint triangles,
$\nu_t(G)$ and $\nu_t^*(G)$ are asymptotically the same. It follows from these results that:
\begin{theo}
\label{t:tuza-dense}
$\tau_t(G) \le 2\nu_t(G) + o(n^2)$.
\end{theo}
In light of the fact that Tuza's conjecture is optimal, it is interesting to ask whether the constant $2$ in Theorem \ref{t:tuza-dense} is also optimal
(notice that this question becomes nontrivial for dense graphs with $\tau_t(G)=\Theta(n^2)$).
Perhaps the most interesting case to consider is when $\tau_t(G)$ is as large as one can expect it to be.

It is well known that every graph with $m$ edges can be made bipartite by removing from it less than $m/2$ edges (see \cite{Al96} for the tightest known bounds). In particular, $\tau_t(G) \le m/2 -o(m)$. On  the other hand, there are many different types of graphs that are hard to make triangle-free,
that is, graphs for which $\tau_t(G) \ge m/2 -o(m)$.
For example, complete graphs are hard to make triangle-free, and (sufficiently dense) random graphs are hard to make triangle-free.
It is also easy to construct many other families of graphs that are hard to make triangle-free.
Let us formalize this notion. We say that a graph $G$ is {\em $(1-\delta)$-hard to make $\Delta$-free} if
$\tau_t(G) \ge (1-\delta)(m/2)$.

The following is an immediate consequence of Theorem \ref{t:tuza-dense}.
\begin{coro}
\label{c:tuza-dense}
Let $G$ be a graph with $m$ edges that is $(1-o_n(1))$-hard to make $\Delta$-free. Then,
$$
\nu_t(G) \ge \frac{m}{4}-o(n^2)\;.
$$
\end{coro}

We conjecture that Corollary \ref{c:tuza-dense} is {\em not} optimal, and that $m/4$ can be replaced with $m/3$.
Formally, we conjecture the following.
\begin{conj}
\label{conj:hard}
For every $\epsilon > 0$ and  $\beta > 0$ there exist $N=N(\epsilon,\beta)$ and $\delta=\delta(\epsilon,\beta)$ such that
for all graph with $n > N$ vertices and with $m \ge \beta n^2$ edges that are $(1-\delta)$-hard to make $\Delta$-free,
$$
\nu_t(G) \ge (1-\epsilon)\frac{m}{3}\;.
$$
\end{conj}
Since for any $m$-edge graph we have $\nu_t(G) \le m/3$, Conjecture \ref{conj:hard} states that dense graphs that are hard to make $\Delta$-free
have an asymptotically optimal triangle packing: all but a negligible fraction of the edges are packed.

A weakened, but still challenging version of Conjecture \ref{conj:hard}, asks for a constant improvement over the $1/4$ bound in Corollary
\ref{c:tuza-dense}.

\begin{conj}
\label{conj:med}
There exists $\alpha > 0$ so that for all $\beta > 0$ there exist $N=N(\beta)$ and $\delta=\delta(\beta)$
such that for all graph with $n > N$ vertices and with $m \ge \beta n^2$ edges that are $(1-\delta)$-hard to make $\Delta$-free,
$$
\nu_t(G) \ge (1+\alpha)\frac{m}{4}\;.
$$
\end{conj}

A further weakening of Conjecture \ref{conj:med} allows the improvement $\alpha$ to depend on the density $\beta$.
The main result of this paper proves that such an improvement always exists.
Hence, for any fixed density, our main result shows that the constant $1/4$ in Corollary \ref{c:tuza-dense} is not optimal, and can be replaced
with a larger constant.
\begin{theo}
\label{t:main}
For every $\beta > 0$ there are $N=N(\beta)$ and $\delta=\delta(\beta)$
such that for all graph with $n > N$ vertices and with $m \ge \beta n^2$ edges that are $(1-\delta)$-hard to make $\Delta$-free,
$$
\nu_t(G) \ge \left(1+\frac{\beta^2}{800}\right)\frac{m}{4}\;.
$$
\end{theo}
The constant $800$ in Theorem \ref{t:main} is by no means optimal and it can be somewhat reduced at the price of complicating the calculations.
Since this has no qualitative impact on the statement of Theorem \ref{t:main}, we make no effort to optimize it.

The next section contains the proof of Theorem \ref{t:main}.
Section 3 considers larger cliques. We prove a bound for the covering number of $K_k$ in terms of the fractional covering number of $K_k$ and
then use it to extend Theorem \ref{t:main} to larger cliques.
Section 4 contains some concluding remarks: a sketch of a generalization of Theorem \ref{t:main} to larger odd cycles,
and an improved integrality gap for the problem of ``maximal triangle-free subgraph'' in dense graphs.

\section{Packing triangles in graphs that are hard to make triangle free}

Since $\nu_t^*(G)=\tau_t^*(G)$ and since, by the result of Haxell and R\"odl mentioned earlier we have $\nu_t^*(G) \le \nu_t(G)+o(n^2)$,
proving Theorem \ref{t:main} is equivalent to proving the following theorem.
\begin{theo}
\label{t:main-2}
For every $\beta > 0$ there are $N=N(\beta)$ and $\delta=\delta(\beta)$
such that for all graph with $n > N$ vertices and with $m \ge \beta n^2$ edges that are $(1-\delta)$-hard to make $\Delta$-free,
$$
\tau_t^*(G) \ge \left(1+\frac{\beta^2}{800}\right)\frac{m}{4}\;.
$$
\end{theo}
We will therefore prove Theorem \ref{t:main-2}, and hence obtain a proof for Theorem \ref{t:main} as well.

We first need to recall some known facts from linear programming. For a graph $G$, let $E(G)$ and $T(G)$ denote the sets of edges and triangles
of $G$, respectively. Let $f: E(G) \rightarrow [0,1]$ be a minimum fractional triangle cover so that
$\sum_{e \in E(G)} f(e) = \tau_t^*(G)$, and let $g: T(G) \rightarrow [0,1]$ be a maximum fractional triangle packing so that
$\sum_{t \in T(G)} g(t) = \nu_t^*(G)$. Then, the duality theorem of linear programming states that $\tau_t^*(G)=\nu_t^*(G)$ and (one of) the complementary slackness conditions states that:
\begin{equation}
\label{e:slackness}
f(e) > 0 {\rm~implies~} \sum_{t \ni e} g(t) = 1\;.
\end{equation}

We designate two sets of edges.
\begin{itemize}
\item
Let $F_0 \subset E(G)$ be $F_0 = \{e~|~f(e)=0\}$.
\item
Let $F_1 \subset E(G)$ be $F_1 = \{e~|~f(e)=1\}$.
\end{itemize}
The proof of Theorem \ref{t:main-2} is split into three cases, according to the cardinalities of $F_0$ and $F_1$.
The first two cases are easy. The first is when $F_1$ is relatively large and the second is when $F_0$ is relatively small.
The remaining case, where $F_1$ is relatively small and $F_0$ is relatively large, is more difficult.
It will be convenient to assume, without loss of generality, that $m=\beta n^2$. Observe that this immediately implies the proof
for $m \ge \beta n^2$.
\vspace{10pt}

\noindent
{\bf Case 1: $|F_1| > \left(\delta+\frac{\beta^2}{800}\right)m/2 $.}

\noindent
Define $G_1 = G(V,E \setminus F_1)$ to be the graph obtained from $G$ by deleting the edges having weight $1$.
We observe that:
\begin{eqnarray}
\label{e1} \tau_t(G_1) & \ge & \tau_t(G)-|F_1|\;. \\
\label{e2} \tau_t^*(G_1) & \le & \tau_t^*(G)-|F_1|\;. \\
\label{e3} \tau_t^*(G_1) & \ge & \frac{1}{2}\tau_t(G_1)\;.
\end{eqnarray}
Indeed, (\ref{e1}) holds since we have deleted $|F_1|$ edges,
(\ref{e2}) holds since the total deleted weight is $|F_1|$,
and (\ref{e3}) holds by Theorem \ref{t:kriv}.
Using these inequalities and the assumption on the size of $F_1$ we have:
\begin{eqnarray*}
\tau_t^*(G) & \ge & \tau_t^*(G_1)+|F_1| \\
& \ge & \frac{1}{2}\tau_t(G_1)+|F_1|\\
& \ge & \frac{1}{2}(\tau_t(G)-|F_1|) + |F_1|\\
& = & \frac{1}{2}\tau_t(G)+ \frac{1}{2}|F_1|\\
& \ge & \frac{1}{2}(1-\delta)\frac{m}{2}+\left(\delta+\frac{\beta^2}{800}\right)\frac{m}{4}\\
& = & \left(1+\frac{\beta^2}{800}\right)\frac{m}{4}\;.
\end{eqnarray*}

\vspace{10pt}

\noindent
{\bf Case 2: $|F_0| < \left(1-\frac{3\beta^2}{800}\right)m/4$.}

\noindent
The complementary slackness condition (\ref{e:slackness}) implies that
$$
\sum_{e \in E \setminus F_0} \sum_{t \ni e} g(t) = m-|F_0|\;.
$$
As each triangle is counted at most three times we have that
$$
\tau_t^*(G)=\nu_t^*(G) \ge \frac{1}{3}(m-|F_0|)\;.
$$
Using the assumption on the size of $F_0$ we obtain:
$$
\tau_t^*(G) \ge \frac{1}{3}\left(m-\left(1-\frac{3\beta^2}{800}\right)m/4\right) = \left(1+\frac{\beta^2}{800}\right)\frac{m}{4}\;.
$$

\vspace{10pt}

\noindent
{\bf Case 3: $|F_0| \ge \left(1-\frac{3\beta^2}{800}\right)m/4$ and $|F_1| \le \left(\delta+\frac{\beta^2}{800}\right)m/2 $.}

\noindent
Choosing $\delta=\delta(\beta)$ to be sufficiently small, and since $\beta \le 1/2$, our assumptions in this case imply in particular that
\begin{eqnarray}
\label{e4} |F_0| & \ge & 0.248m\;. \\
\label{e5} |F_1| & \le & \beta^2 m/1599\;.
\end{eqnarray}

Consider the graph $H=G(V,F_0)$ consisting only of the edges having weight zero.
Notice that $H$ is still dense as it has at least $0.248m$ edges, and that $H$ is triangle-free since otherwise $f$ would not have been
a fractional triangle cover. The following lemma proves (and quantifies) that $H$ has a dense {\em induced} bipartite subgraph, where each vertex class
is a (partial) neighborhood.
\begin{lemma}
\label{l:bip}
$H$ has an induced bipartite subgraph $(A \cup B, F^*)$ with $|F^*| \ge \beta^2 m/500$. Furthermore $A$ has a common neighbor in $H$ and $B$ has a common neighbor in $H$.
\end{lemma}
{\bf Proof:}\,
Since $H$ has at least $0.248m$ edges, if we repeatedly delete vertices with degree less than $0.124\beta n$ we remain with a
subgraph $H'$ with at least $0.124m$ edges and minimum degree at least $0.124\beta n$.
We will prove that the lemma holds already for $A$ and $B$ in $H'$.

Consider a random set $C$ of $c=\lceil 1/(0.124\beta) \rceil$ vertices, chosen from the vertex set $V$.
We say that a vertex $x$ of $H'$ is dominated by $C$ if some vertex of $C$ is a neighbor of $x$ in $H'$.
Clearly, the probability that a vertex $x \in H$ is {\em not} dominated by $C$ is less than $(1-c/n)^{d_x}$.
Let $Q$ denote the set of all edges of $H'$ that are incident with vertices that are not dominated by $C$.
Hence, the expected size of $Q$ satisfies:
$$
E[|Q|] < \sum_{x \in H'} d_x\left(1-\frac{c}{n}\right)^{d_x}\;.
$$
Recall that $d_x \ge 0.124\beta n$ and notice that a term in the last inequality is maximized for $d_x = 0.124\beta n$.
Thus,
$$
E[|Q|] < \sum_{x \in H'} 0.124\beta n\left(1-\frac{c}{n}\right)^{0.124\beta n} \le \sum_{x \in H'} \frac{0.124\beta n}{e} < \frac{0.124\beta n^2}{e}\;.
$$
In particular, there exists a choice of $C$ so that after removing from $H'$ the vertices that are not dominated by $C$ we remain with a subgraph
$H''$ whose number of edges is at least
$$
0.124m- \frac{0.124\beta n^2}{e} = 0.124m- \frac{0.124m}{e} =0.124(1-1/e)m\;.
$$

As each vertex of $H''$ is dominated by $C$, let us select for each vertex of $H''$, a vertex of $C$ that dominates it.
This partitions the vertices of $H''$ into $c$ parts $\{A_u ~|~ u \in C\}$ where $A_u$ consists of the vertices of $H''$ that chose
$u$ as their dominating vertex. Observe that $A_u$ induces an independent set in $H$, since otherwise we would have, together with $u$,
a triangle in $H$, contradicting its triangle-freeness. Consider now all the pairs $(A_u,A_w)$ for $u,w \in C$ and $u \neq w$.
Each edge of $H''$ appears in precisely one of these pairs. Since $H''$ has more than $0.124(1-1/e)m$ edges and since there are only
${c \choose 2}$ pairs, we must have a pair $(A_u,A_w)$ so that the bipartite subgraph induced by it contains at least
$$
\frac{0.124(1-1/e)m}{{c \choose 2}} > \frac{m\beta^2}{500}
$$
edges. Since $A_u$ has $u$ as a common neighbor and since $A_w$ has $w$ as a common neighbor, the lemma follows. \qed

Let us consider the subgraph induced by $A$ and $B$ of Lemma \ref{l:bip} in $G$.
By the lemma we have that $|E(A,B)| \ge m\beta^2/500$.
We claim that $E(A)$ and $E(B)$ (edges with both endpoints in $A$ or both endpoints in $B$) contain only edges of $F_1$.
Indeed, recall that by Lemma \ref{l:bip}, $A$ has a common neighbor, say, $v$, with all of the edges $(v,a)$ with $a \in A$ being in
$F_0$. Now, if we had an edge $(a,a') \in E(A)$ with $f((a,a')) < 1$ then the triangle $(a,a',v)$ would have total weight less than $1$,
contradicting the fact that $f$ is a fractional triangle cover. An identical argument holds for $E(B)$.
It follows that
$$
|E(A,B)| - |E(A)| - |E(B)| \ge \frac{m\beta^2}{500} - |F_1| \ge \frac{m\beta^2}{800}\;.
$$

To complete the proof of Case 3, we proceed as follows. We split the vertices of $V \setminus (A \cup B)$ into two parts $X$ and $Y$
at random. We consider the cut $(A \cup X~,~B \cup Y)$ and compute the expected number of edges crossing it.
Each edge of $E(A,B)$ crosses it by definition. On the other hand, each edge with at least one endpoint in $X \cup Y$ crosses it
with probability $1/2$. Hence the expected size of this cut is
$$
|E(A,B)| + \frac{1}{2}\left(m-|E(A,B)|-|E(A)|-|E(B)|\right) \ge \frac{m}{2} + \frac{m\beta^2}{1600}\;.
$$
Hence, such a cut exists, implying that we can remove from $G$ the non-edges of this cut to obtain a triangle-free (in fact, bipartite)
subgraph. The number of edges thus removed is less than
$$
m - \frac{m}{2} - \frac{m\beta^2}{1600} = \frac{m}{2}(1 - \frac{\beta^2}{800}) < \frac{m}{2}(1-\delta)
$$
contradicting the assumption that $\tau_t(G) \ge \frac{m}{2}(1-\delta)$.
This completes the proof of Theorem \ref{t:main-2}, which, as noted earlier, implies Theorem \ref{t:main}. \qed

\section{Larger cliques}

Throughout this section we fix $k \ge 4$, we let
$\tau_k(G)$ denote the minimum size of a $K_k$-cover, and let $\nu_k(G)$ denote the maximum size of a $K_k$-packing of a graph $G$.
The trivial bounds in this case are $\nu_k(G) \le \tau_k(G) \le {k \choose 2}\nu_k(G)$.

Denoting by $\tau_k^*(G)$ and $\nu_k^*(G)$ the respective (and equal) fractional parameters, Krivelevich's proof for triangles \cite{Kr95} can be generalized to yield:
\begin{equation}
\label{e:kk-kriv}
\tau_k(G) \le \left({k \choose 2}-1\right)\tau_k^*(G)\;.
\end{equation}
We omit the details of this easy generalization since the bounds we shall obtain in this section are better.

As for the case of triangles, the theorem of Haxell and R\"odl \cite{HaRo01} asserts that $|\nu_k(G)-\nu^*_k(G)| = o(n^2)$.
Thus, an immediate corollary analogous to Theorem \ref{t:tuza-dense} is:
\begin{coro}
\label{c:kk-dense}
$\tau_k(G) \le \left({k \choose 2}-1\right)\nu_k(G) + o(n^2)$.
\end{coro}

The goal of this section is to prove a significantly better bound, replacing ${k \choose 2}-1$ with a much smaller value.
We shall do that by improving upon (\ref{e:kk-kriv}).

\begin{theo}
\label{t:kk}
$$\tau_k(G) \le \lfloor k^2/4 \rfloor\tau_k^*(G)\;.$$
\end{theo}
{\bf Proof:}\,
Consider the following process which creates a sequence of spanning subgraphs $G_i$ of $G$, starting with $G=G_0$.
Each $G_i$ is obtained from its predecessor $G_{i-1}$ by deleting a single edge according to the rule specified below.
We will halt this process one this rule cannot be applied.
We denote the final graph in our sequence by $G_t$. Hence we have $0 \le t \le m$.

Let $f_i$ and $g_i$ be a minimum fractional $K_k$-cover and a maximum fractional $K_k$-packing of $G_i$, respectively.
Assume first that some $K_k$ of $G_i$ contains ${k \choose 2}- \lfloor k^2/4 \rfloor$ edges that are assigned weight $0$ by $f_i$.
This means that the total weight of the remaining $\lfloor k^2/4 \rfloor$ edges of this $K_k$ is at least $1$, so there is some
edge $e_i$ with $f_i(e_i) \ge 1/\lfloor k^2/4 \rfloor$. We let $G_{i+1}= G_i - \{e_i\}$.
If no $K_k$ of $G_i$ contains ${k \choose 2}- \lfloor k^2/4 \rfloor$ edges that are assigned weight $0$ by $f_i$ then we halt the
sequence and $G_i=G_t$ is the final graph in the sequence.

We observe the following inequalities:
\begin{eqnarray}
\label{e6} \tau_k(G_t) & \ge & \tau_k(G)- t\;. \\
\label{e7} \tau_k^*(G_t) & \le & \tau_k^*(G)-\frac{t}{\lfloor k^2/4 \rfloor}\;. \\
\label{e8} \tau_k^*(G_t) & \ge & \frac{\tau_k(G_t)}{{k \choose 2}-1}\;.
\end{eqnarray}
Indeed, (\ref{e6}) holds since we have deleted $t$ edges to get from $G$ to $G_t$,
(\ref{e7}) holds since $\tau_k^*(G_{i+1}) \le \tau_k^*(G_i)-1/\lfloor k^2/4 \rfloor$,
and (\ref{e8}) holds by (\ref{e:kk-kriv}).
Using these inequalities we have:
\begin{eqnarray}
\label{e9}
\nonumber \tau_k^*(G) & \ge & \tau_k^*(G_t)+ \frac{t}{\lfloor k^2/4 \rfloor}\\
\nonumber & \ge & \frac{\tau_k(G_t)}{{k \choose 2}-1}+\frac{t}{\lfloor k^2/4 \rfloor}\\
\nonumber & \ge & \frac{\tau_k(G)-t}{{k \choose 2}-1} + \frac{t}{\lfloor k^2/4 \rfloor}\\
& = & \frac{\tau_k(G)}{{k \choose 2}-1} -\frac{t}{{k \choose 2}-1} + \frac{t}{\lfloor k^2/4 \rfloor}\;.
\end{eqnarray}

Let $0 \le \alpha \le 1$  be a parameter.
Consider first the case where $G_t$ has at most $\alpha(m-t)$ edges that are assigned weight $0$ by $f_i$.
Thus, at least $(1-\alpha)(m-t)$ edges of $G_t$ are assigned positive weight and, using complementary slackness as in (\ref{e:slackness}) we obtain:
$$
\sum_{e:f_i(e) > 0} \sum_{H \in e}g_i(H) = \sum_{e:f_i(e) > 0} 1 \ge (1-\alpha)(m-t)\;.
$$
(Here the internal sum ranges over all graph $H$ in $G_t$ that are isomorphic to $K_k$.)
Since $K_k$ has ${k \choose 2}$ edges, this implies, in particular:
$$
\tau_k^*(G_t)=\nu_k^*(G_t) \ge \frac{(1-\alpha)(m-t)}{{k \choose 2}}\;.
$$
By (\ref{e7}) we have:
\begin{equation}
\label{e10}
\tau_k^*(G) \ge \frac{(1-\alpha)(m-t)}{{k \choose 2}}+\frac{t}{\lfloor k^2/4 \rfloor}\;.
\end{equation}

Consider next the case where $G_t$ has at least $\alpha(m-t)$ edges that are assigned weight $0$ by $f_i$.
This means that the spanning subgraph $P$ of $G_t$ consisting of the edges having positive weight has
at most $(1-\alpha)(m-t)$ edges. Since any graph can be made bipartite by removing less than half of its edges,
we can delete from $P$ a subset $F$ of less than $(1-\alpha)(m-t)/2$ edges to make $P$ bipartite.

We claim that the spanning subgraph $Q$ of $G_t$ obtained by removing $F$ from $G_t$ is $K_k$-free.
Assume that $Q$ has a $K_k$. The edges with positive weight form a bipartite subgraph on $k$ vertices inside
this $K_k$. By Mantel's Theorem \cite{Ma07}, the number of such edges is at most $\lfloor k^2/4 \rfloor$.
This implies that this $K_k$ contains at least ${k \choose 2}-\lfloor k^2/4 \rfloor$ edges with zero weight,
contradicting the fact that $G_t$ was the last graph in the sequence and has no copy of $K_k$
with this amount of zero weight edges.
We have therefore proved:
$$
\tau_k(G_t) \le \frac{(1-\alpha)(m-t)}{2}\;.
$$
By (\ref{e6}) we have:
\begin{equation}
\label{e11}
\tau_k(G) \le \frac{(1-\alpha)(m-t)}{2} + t\;.
\end{equation}
By (\ref{e10}) and (\ref{e11}) we have:
\begin{equation}
\label{e12}
\tau_k^*(G) \ge \frac{2\tau_k(G)-2t}{{k \choose 2}}+\frac{t}{\lfloor k^2/4 \rfloor}\;.
\end{equation}

So, (\ref{e9}) and (\ref{e12}) both supply lower bounds for $\tau_k^*(G)$ in terms of $\tau_k(G)$ and $t$.
In particular, the maximum of both bounds can be used as a lower bound for $\tau_k^*(G)$.
For $k \ge 4$ observe that (\ref{e9}) increases as $t$ increases and (\ref{e12}) decreases as $t$ increases.
Hence, the maximum of both bounds is minimized when they are equal which, in turn,
happens when $t=\tau_k(G)$. In this extremal point we have
$$
\tau_k^*(G) \ge \frac{\tau_k(G)}{\lfloor k^2/4 \rfloor}\;.
$$
Thus, $\tau_k(G) \le \lfloor k^2/4 \rfloor\tau_k^*(G)$ proving the theorem. \qed

Theorem \ref{t:kk} immediately implies the following improvement of Corollary \ref{c:kk-dense}
\begin{coro}
\label{c:kk-dense-2}
$\tau_k(G) \le \lfloor k^2/4 \rfloor\nu_k(G) + o(n^2)$.
\end{coro}

It is well known that every graph with $m$ edges can be made $(k-1)$-partite by removing from it less than $m/(k-1)$ edges.
One just considers a random partition of the vertex set into $k-1$ parts and observes that the probability of an edge
having both of its endpoints in the same part is less than $1/(k-1)$.
In particular, $\tau_k(G) \le m/(k-1) -o(m)$. As for the case of triangles, there are many different types of graphs
that are hard to make $K_k$-free, that is, graphs for which $\tau_k(G) \ge m/(k-1) -o(m)$.
We thus say that a graph $G$ is {\em $(1-\delta)$-hard to make $K_k$-free} if
$\tau_k(G) \ge (1-\delta)m/(k-1)$.

The following is an immediate consequence of Corollary \ref{c:kk-dense-2}.
\begin{coro}
\label{c:kk-hard}
Let $G$ be a graph with $m$ edges that is $(1-o_n(1))$-hard to make $K_k$-free. Then,
$$
\nu_k(G) \ge \frac{m}{(k-1)\lfloor k^2/4 \rfloor}-o(n^2)\;.
$$
\end{coro}
Observe that for, say, $K_4$ we get that dense graphs that are hard to make $K_4$-free have roughly $m/12$ edge-disjoint copies
of $K_4$. As each $K_4$ has $6$ edges this implies that a fraction of roughly $1/2$ of the edges can be packed with edge-disjoint copies of $K_4$.
More generally, for $K_k$, we get that a fraction of $2/k$ of the edges can be packed with edge-disjoint copies of $K_k$ (if $k$ is odd then this
fraction is a bit larger). It is plausible that conjecture \ref{conj:hard} can be extended from triangles to larger cliques.
That is, all but a negligible fraction of the edges can be packed with edge-disjoint copies of $K_k$. 

\section{Concluding remarks}

The proof of Theorem \ref{t:main} can be extended to other odd cycles.
Denoting the covering and packing numbers by $\tau_{C_k}(G)$ and $\nu_{C_k}(G)$ respectively, the analogous result states that
for $\beta$-dense graphs that are $(1-\delta)$-hard to make $C_k$-free one has
$$
\nu_{C_k}(G) \ge \left(1+c\beta^2\right)\frac{m}{2k-2}\;.
$$
where $c$ is an absolute constant.
Observe that for any graph $G$ we have $\nu_{C_k}(G) \le m/k$.

The proof is essentially the same with the following minor differences.
We use a straightforward extension of the result of Krivelevich for cycles of length $k$, which
states that $\tau_{C_k}(G) \le (k-1)\tau^*_{C_k}(G)$ (see also \cite{KoLaNu08} for this observation),
and the result of Haxell and R\"odl applied to $C_k$ stating that $|\nu_{C_k}(G)-\nu^*_{C_k}(G)| = o(n^2)$.
As in the proof of Theorem \ref{t:main-2}, we split into three cases according to the relative sizes of $F_0$ and $F_{1/(k-2)}$
where the latter are all edges with weight at least $1/(k-2)$. Observe that this coincides with the definition of $F_1$ for the case of triangles.
The only real difference is in Case 3. In Lemma \ref{l:bip} we can no longer claim that $H$ is triangle-free. Rather, it is $C_k$-free.
This means that any neighborhood of a vertex is no longer forced to be an independent set, but rather it is forced not to
contain a path of length $k-2$. But this, in turn, implies that each neighborhood in $H$ is sparse and has only a linear number
of edges, which is negligible in the dense setting. Also, when using Lemma \ref{l:bip} by looking at the subgraph induced by $A \cup B$ in $G$,
we can no longer claim that it contains only edges of $F_{1/(k-2)}$ with both endpoints in $A$ or both endpoints in $B$.
However, it certainly does not contain a path of length $k-2$ of edges not in $F_{1/(k-2)}$ with both endpoints in the same class.
Thus, there are only a negligible (linear) number of edges not in $F_{1/(k-2)}$ that are inside $A$ or inside $B$.
Hence, the same argument as in Case 3 for triangles, also holds here.

\vspace{10pt}
Theorem \ref{t:kk} supplies, in particular, an efficient approximation algorithm for the NP-Hard problem of
computing $\tau_k(G)$. Its approximation ratio is $\lfloor k^2/4 \rfloor$. It also bounds the integrality gap of
this problem by $\lfloor k^2/4 \rfloor$.

\vspace{10pt}
Consider the problem of finding a maximal triangle-free subgraph. Its fractional relaxation is thus to
assign weights in $[0,1]$ to the edges so that for each triangle, the sum of the weights is not larger than $2$.
The goal is to maximize the sum of the weights of such an assignment.
Denoting the corresponding parameters by $\rho_t(G)$ and $\rho_t^*(G)$ we have, by definition, $\rho_t(G) = m - \tau_t(G)$ and
$\rho_t^*(G) = m - \tau_t^*(G)$.
The (asymptotic) integrality gap of this problem is known to be between $1.5$ and $4/3$. The lower bound comes from the complete graph:
The integral solution is $n^2/4(1-o(1))$, while the fractional solution comes from assigning a weight of $2/3$ to each edge,
thereby obtaining total weight of $n^2/3(1-o(1))$. The upper bound follows from the Krivelevich's result $\tau_t(G) \le 2\tau_t^*(G)$ after
some easy arithmetic manipulations.

Our proof of Theorem \ref{t:main-2} improves upon the upper bound for dense graphs.
Suppose that $G$ is a graph with $m=\beta n^2$ edges. Assume first that $\tau_t(G) \ge (1-\delta)m/2$.
By Theorem \ref{t:main-2}, $\rho_t^*(G) \le 3m/4 - m\beta^2/3200$. On the other hand, for any graph we have $\rho_t(G) \ge m/2$.
Thus, the integrality gap in this case is at most $3/2-\beta^2/1600$.
Consider next the case $\tau_t(G) \le (1-\delta)m/2$. Hence, $\rho_t(G) \ge (1+\delta)m/2$.
By Krivelevich's result, we have
$$
\rho_t(G) = m - \tau_t(G) \ge m - 2\tau_t^*(G) = 2\rho_t^*(G) -m\;.
$$
This implies that the integrality gap is at most $1/2+m/(2\rho_t(G))$.
In our case this implies an integrality gap of $1/2+1/(1+\delta)=3/2-\delta/(1+\delta)$.
Observe that from the proof the theorem, we actually have that the function $\delta(\beta)$ in the statement of the theorem can be chosen to be $\delta=c\beta^2$ where $c$ is a small absolute constant.

\section*{Acknowledgment}

I thank Guy Kortsarz and Michael Krivelevich for useful discussions.

\end{document}